\documentclass[12pt]{article}
\author{Gun Srijuntongsiri\thanks{Sirindhorn International Institute of Technology, 131 Moo 5, Tiwanont Road,
Bangkadi, Muang, Pathumthani, 12000, Thailand. Email:
gun@siit.tu.ac.th.},
Stephen A. Vavasis\thanks{MC 6054, University of Waterloo, 200
University Avenue W., Waterloo, ON N2L 3G1, Canada. Email:
vavasis@math.uwaterloo.ca.}}
\title{Properties of polynomial bases used in a line-surface intersection algorithm\thanks{Supported in part by NSF DMS 0434338 and NSF CCF 0085969.}}
\date{\today}

\usepackage{amsmath, amsthm, amssymb,graphicx,subfig}

\begin{document}
\maketitle
\newcommand{\cond}[1]{\mathop{\rm{cond}}(#1)}
\newcommand{\norm}[1]{\left\|#1\right\|}
\newcommand{\normtwo}[1]{\left\|#1\right\|_2}
\newcommand{\abs}[1]{\left|#1\right|} 
\newcommand{\norma}[1]{\left\|#1\right\|}   
\newcommand{\inv}[1]{#1^{-1}}
\newcommand{\nchoosek}[2]{\left(\begin{array}{c} #1\\ #2 \end{array} \right)}

\newtheorem{thm}{Theorem}[section]
\newtheorem{cor}[thm]{Corollary}
\newtheorem{lemma}[thm]{Lemma}
\newtheorem{prop}{Proposition}[thm]

\theoremstyle{remark}
\newtheorem{rem}[thm]{Remark}

\begin{abstract}
In \cite{srijuntongsiri_lsi}, Srijuntongsiri and Vavasis propose
the \emph{Kantorovich-Test Subdivision algorithm}, or KTS, which is
an algorithm for finding all zeros of a polynomial system in a bounded region
of the plane.  This algorithm can be used to find the intersections between a
line and a surface.  The main features of KTS are that it can operate on
polynomials represented in any basis that satisfies certain
conditions and that its efficiency has an upper bound
that depends only on the conditioning of the problem and the choice of
the basis representing the polynomial system.

This article explores in detail the dependence of the efficiency of the KTS algorithm on the choice of basis.
Three bases are considered: the power, the Bernstein, and the Chebyshev bases.
These three bases satisfy the basis properties required by KTS.
Theoretically, Chebyshev case has the smallest upper bound on its running time.
The computational results, however, do not show that Chebyshev case performs
better than the other two.
\end{abstract}

\section{The line-surface intersection problem and the required basis properties}
\label{section_basis_prop}

Let $\phi_0$, \ldots, $\phi_n$ denote a basis for the set of
univariate polynomials of degree at most $n$.  For example,
the power basis is defined by $\phi_i(t)=t^i$.  The line-surface
intersection problem can be reduced to the problem
of finding all zeros of
\begin{equation}
\label{maineq}
\begin{array}{llll}
f(u,v) & \equiv & \displaystyle\sum_{i=0}^m \sum_{j=0}^n c_{ij}
\phi_i(u)\phi_j(v), & 0 \leq u,v \leq 1,
\end{array}
\end{equation}
where $c_{ij} \in \mathbb{R}^2$ $(i=0,1,\ldots,m;
j=0,1,\ldots,n)$ denote the coefficients \cite{srijuntongsiri_lsi}.

For this article, let the
notation $\norm{\cdot}$ refer specifically to infinity
norm.  Other norms are explicitly notated so.
The \emph{Kantorovich-Test Subdivision algorithm} (KTS
in short), proposed by Srijuntongsiri and Vavasis \cite{srijuntongsiri_lsi},
can be used to solve (\ref{maineq}).
KTS works with any polynomial basis $\phi_i(u)\phi_j(v)$
provided that the following properties hold:
\begin{enumerate}
\item There is a natural interval $[l,h]$ that is the domain for the polynomial.
In the case of Bernstein polynomials, this is $[0,1]$, and in the case of power and
Chebyshev polynomials, this is $[-1,1]$.

\item It is possible to compute a bounding polytope $P$ of $S =
\{f(u,v) : l \leq u,v \leq h \}$, where $f(u,v) = \sum_{i=0}^m
\sum_{j=0}^n c_{ij} \phi_i(u)\phi_j(v)$ and $c_{ij} \in
\mathbb{R}^d$ for any $d \geq 1$, that satisfies the following
properties: \label{bounding_prop}
\begin{enumerate}
\item Determining whether $0 \in P$ can be done efficiently
(ideally in $O(mn)$ operations).

\item The polytope $P$ is affinely invariant. In other words,
the bounding polytope of $\{Af(u,v)+b: l \leq u,v \leq h\}$
 is $\{Ax+b : x \in P \}$ for any nonsingular matrix
$A \in \mathbb{R}^{d \times d}$ and any vector $b \in \mathbb{R}^d$.

\item For any $y \in P$,
\begin{equation}
\norm{y} \leq \theta\max_{l \leq u,v \leq h} \norm{ f(u,v) },
\label{thetadef}
\end{equation}
where $\theta$ is a function of $m$ and $n$.
\label{prop_theta}

\item If $d=1$, then the endpoints of $P$ can be computed efficiently (ideally in $O(mn)$ time).
\label{prop_endpoint}
\end{enumerate}

\item It is possible to reparametrize with $[l,h]^2$ the surface
$S_1 = \{ f(x) : x \in \bar{B}(x^0,r) \}$, where $x^0 \in \mathbb{R}^2$ and $r \in \mathbb{R} > 0$.
In other words, it is possible (and efficient) to compute the polynomial
$\hat{f}$ represented in the same basis such that $S_1=\{\hat f(\hat x): \hat x \in [l,h]^2\}$.

\item Constant polynomials are easy to represent. \label{constant_poly_easy}

\item Derivatives of polynomials are easy to determine in the same basis. (preferably in $O(mn)$ operations).
\label{prop_deriv}
\end{enumerate}
We are generally interested in the case where $d=2$.  In this case, we call $P$ a bounding polygon.  Recall that $P$
is a bounding polygon of $S$ if and only if $x \in S$ implies $x
\in P$.

\section{The Kantorovich-Test Subdivision algorithm}

The description of KTS, as well as the definitions of
the quantities mentioned in the description, are given below.
More details can be found in \cite{srijuntongsiri_lsi}.

For a given zero $x^*$ of polynomial $f$,
let $\omega_*(x^*)$ and $\rho_*(x^*)$
be quantities satisfying the conditions that, first, $\omega_*(x^*)$ is the smallest Lipschitz constant
for $f'(x^*)^{-1}f'$, i.e.,
\begin{equation}
\label{lip1}
\norm{f'(x^*)^{-1}\left( f'(x)-f'(y) \right)} \leq \omega_*(x^*) \cdot \norm{x-y} \textrm{ for all } x,y \in \bar{B}(x^*,\rho_*(x^*))
\end{equation}
and, second,
\begin{equation}
\label{rhostar}
\rho_*(x^*) = \frac{2}{\omega_*(x^*)}.
\end{equation}
Define
\[
\gamma(\theta) = 1 / \left(4\sqrt{\theta(4\theta+1)}-8\theta \right),
\]
where $\theta$ is as in (\ref{thetadef}).
Define $\omega_{D'}$ to be
the smallest nonnegative constant $\omega$ satisfying
\begin{equation}
\begin{array}{llll}
\norm{f'(x^*)^{-1}\left(f'(y)-f'(z)\right)} & \leq & \omega \cdot \norm{y-z}, & y,z \in D', x^* \in [0,1]^2
 \\
& & & \textrm{satisfying } f(x^*) = 0,  \label{om}\end{array}
\end{equation}
where
\begin{equation}
\label{dprimedef} D' =
\left[-\gamma(\theta),1+\gamma(\theta)\right]^2.
\end{equation}
Denote $\omega_f$ as the maximum of
$\omega_{D'}$ and all $\omega_*(x^*)$
\[
\omega_f = \max\{\omega_{D'}, \max_{x^* \in \mathbb{C}^2 : f(x^*)=0} \omega_*(x^*)\}.
\]
Finally, define the condition number of $f$ to be
\begin{equation}
\label{conddef}
\cond{f} = \max \{ \omega_f, \max_{x^* \in \mathbb{C}^2 : f(x^*)=0,y \in [0,1]^2} \norm{\inv{f'(x^*)}f'(y)}\}.
\end{equation}

We define the \emph{Kantorovich test} on a region $X =
\bar{B}(x^0,r)$ as the application of Kantorovich's Theorem on
the point $x^0$ using $\bar{B}(x^0,2\gamma(\theta)r)$ as the
domain (refer to \cite{deuflhard, kantorovich} for the statement of
Kantorovich's Theorem).  The region $X$ passes the Kantorovich
test if $\eta \omega \leq 1/4$ and $\bar{B}(x^0,\rho_-)
\subseteq D'$.

The other test KTS uses is the exclusion test. For a
given region $X$, let $\hat{f}_X$ be the polynomial in the basis
$\phi_i(u)\phi_j(v)$ that reparametrizes with $[l,h]^2$ the
surface defined by $f$ over $X$.  The region $X$ passes the
\emph{exclusion test} if the bounding polygon of $\{
\hat{f}_X(u,v) : l \leq u,v \leq h \}$ excludes the origin.

Having defined the above prerequisites, the description of KTS
can now be given.

\begin{flushleft}
\textbf{Algorithm KTS}:
\end{flushleft}
\begin{itemize}
\item Let $Q$ be a queue with $[0,1]^2$ as its only entry. Set $S = \emptyset$.
\item Repeat until $Q = \emptyset$
\begin{enumerate}
\item Let $X$ be the patch at the front of $Q$.  Remove $X$ from $Q$.
\item If $X \not\subseteq X_S$ for all $X_S \in S$,
\begin{itemize}
\item Perform the exclusion test on $X=\bar{B}(x^0,r)$
\item If $X$ fails the exclusion test,
\begin{enumerate}
\item Perform the Kantorovich test on $X$
\item If $X$ passes the Kantorovich test,
\begin{enumerate}
\item Perform Newton's method starting from $x^0$ to find a zero $x^*$.
\item If $x^* \not\in X_S$ for any $X_S \in S$ (i.e., $x^*$ has not been found previously),
\begin{itemize}
\item Compute $\rho_*(x^*)$ and its associated $\omega_*(x^*)$ by binary search.
\item Set $S = S \cup \left\{\bar{B}(x^*,\rho_*(x^*))\right\}$.
\end{itemize}
\end{enumerate}
\item Subdivide $X$ along both $u$ and $v$-axes into four equal subregions. Add these subregions
to the end of $Q$.
\end{enumerate}
\end{itemize}
\end{enumerate}
\end{itemize}

The following theorem shows that the efficiency of KTS has an upper bound
that depends only on the conditioning of the problem and the choice of
the basis.
\begin{thm}
\label{thm3} Let $f(x)=f(u,v)$ be a polynomial system in basis
$\phi_i(u)\phi_j(v)$ in two dimensions with generic coefficients
whose zeros are sought.  Let $X = \bar{B}(x^0,r)$ be a patch under
consideration during the course of the KTS algorithm. The
algorithm does not need to subdivide $X$ if
\begin{equation}
\label{thm3ass} r \leq \frac{1}{2}\cdot\min
\left\{\frac{1-1/\gamma(\theta)}{\omega_{D'}},
\frac{1}{2\theta\cond{f}^2}\right\}.
\end{equation}
\begin{proof} See \cite{srijuntongsiri_lsi}.
\end{proof}
\end{thm}
\begin{flushleft}
\textbf{Remark:} Both terms in the bound on the right-hand side of (\ref{thm3ass}) are
increasing as a function of $1/\theta$.
Therefore, our theorem predicts that the KTS algorithm will be more efficient
for $\theta$ as small as possible (close to 1).
\end{flushleft}

\section{Properties of the power, Bernstein, and Chebyshev bases}

As mentioned above, the basis used to represent the polynomial
system must satisfy the properties listed in Section
\ref{section_basis_prop} for KTS to work efficiently. Three bases,
the power, Bernstein, and Chebyshev bases are examined in detail.
The power basis for polynomials of degree $n$ is $\phi_k(t) = t^k$
$(0 \leq k \leq n)$. The Bernstein basis is $\phi_k(t) =
Z_{k,n}(t) = \nchoosek{n}{k}(1-t)^{n-k}t^k$ $(0 \leq k \leq n)$.
The Chebyshev basis is $\phi_k(t) = T_k(t)$  $(0 \leq k \leq n)$,
where $T_k(t)$ is the Chebyshev polynomial of the first kind
generated by the recurrence relation
\begin{eqnarray}
T_0(t) &= &1, \nonumber \\
T_1(t) &= & t, \nonumber \\
T_{k+1}(t) & = & 2tT_k(t)-T_{k-1}(t) \textrm{ for } k \geq 1.
\label{cheb_recurrence}
\end{eqnarray}
Another way to define the Chebyshev polynomials of the first kind
is through the identity
\begin{equation}
\label{cheb_cos_form}
T_k(\cos \alpha) = \cos k\alpha.
\end{equation}
This second definition shows, in particular, that all zeros of $T_k(t)$
lies in $[-1,1]$.  It also shows that $-1 \leq T_k(t) \leq 1$ for any $-1 \leq t \leq 1$.

The rest of this article shows that the power, Bernstein, and Chebyshev
bases all satisfy these basis properties.  The values $\theta$'s of
the three bases and their corresponding bounding polygons are also
derived as these values dictate the efficiency
of KTS operating on such bases.  The upper bound of the efficiency of KTS
is lowest when it operates on the basis with the smallest $\theta$.

\subsection{Bounding polygons}

\label{section_bounding_polygon} The choices of $l$
and $h$ and the definitions of bounding polygons of the surface
$S = \{ f(u,v) : l \leq u,v \leq h\}$, where $f(u,v)$ is represented
by one of the three bases, that satisfy the required properties
are as follows: For Bernstein
basis, the convex hull of the coefficients (control points), call
it $P_1$, satisfies the requirements for $l = 0$ and $h=1$.  The
convex hull $P_1$ can be described as
\begin{equation}
\label{p1def} P_1 = \left\{ \sum_{i,j} c_{ij} s_{ij} : \sum_{i,j}
s_{ij} = 1, 0 \leq s_{ij} \leq 1 \right\}.
\end{equation}
For power and Chebyshev bases, the bounding polygon
\begin{equation}
\label{p2def} P_2 = \left\{ c_{00} + \sum_{i+j > 0} c_{ij} s_{ij} : -1
\leq s_{ij} \leq 1 \right\}
\end{equation}
satisfies the requirements for $l = -1$ and $h = 1$.  Note that
$P_2$ is a bounding polygon of $S$ in the Chebyshev case since $|T_k(t)| \leq 1$
for any $k \geq 0$ and any $t \in [-1,1]$.
Determining whether $0 \in P_2$ is done by solving a small linear programming problem.
To determine if $0 \in P_1$,
the convex hull is constructed by conventional method and is tested if
it contains the origin.

The affine and translational invariance of $P_1$ and $P_2$
for their respective bases can be verified as follows:  Let
\[
g(u,v) = Af(u,v) + b =
\sum_{i=0}^m \sum_{j=0}^n c'_{ij}\phi_i(u)\phi_j(v).
\]
For the Bernstein basis,
by using the property that $\sum_{k=0}^n Z_{k,n}(t) = 1$, it is seen that
$c'_{ij} = Ac_{ij}+b$ for all $c_{ij}$'s.  Therefore, the bounding polygon of
$\{g(u,v): 0 \leq u,v \leq 1 \}$ is
\begin{eqnarray}
P_1' & = & \left\{ \sum_{i,j} c'_{ij} s_{ij} : \sum_{i,j}
s_{ij} = 1, 0 \leq s_{ij} \leq 1 \right\} \nonumber \\
& = & \left\{ \sum_{i,j} (Ac_{ij}+b) s_{ij} : \sum_{i,j}
s_{ij} = 1, 0 \leq s_{ij} \leq 1 \right\} \nonumber \\
& = & \left\{ A\sum_{i,j} c_{ij}s_{ij}+ b\sum_{i,j}s_{ij} : \sum_{i,j}
s_{ij} = 1, 0 \leq s_{ij} \leq 1 \right\} \nonumber \\
& = & \left\{ A\sum_{i,j} c_{ij}s_{ij}+ b : \sum_{i,j}
s_{ij} = 1, 0 \leq s_{ij} \leq 1 \right\} \nonumber \\
& = & \left\{ Ax + b : x \in P_1 \right\}. \nonumber
\end{eqnarray}

For the power and the Chebyshev bases, note that $\phi_0(u)\phi_0(v) = 1$
for both bases.  Hence, $c'_{00} = Ac_{00}+b$ and $c'_{ij} = Ac_{ij}$ for
$i+j > 0$.  The bounding polygon of $\{g(u,v): 0 \leq u,v \leq 1 \}$ for this case is
\begin{eqnarray}
P_2' & = & \left\{ c'_{00} + \sum_{i+j > 0} c'_{ij} s_{ij} : -1
\leq s_{ij} \leq 1 \right\} \nonumber \\
& = & \left\{ Ac_{00}+b + \sum_{i+j > 0} Ac_{ij} s_{ij} : -1
\leq s_{ij} \leq 1 \right\} \nonumber \\
& = & \left\{ A\left( c_{00} + \sum_{i+j > 0} c_{ij} s_{ij} \right) + b: -1
\leq s_{ij} \leq 1 \right\} \nonumber \\
& = & \left\{ Ax + b : x \in P_2 \right\}. \nonumber
\end{eqnarray}

\subsection{The size of the bounding polygons compared to the size of the bounded surface}

Item \ref{prop_theta} of the basis properties in effect ensures that the bounding polygons
are not unboundedly larger than the actual surface itself lest the bounding polygons
lose their usefulness.  The value $\theta$ also can be used as a measure of the
tightness of the bounding polygon.  Recall from Theorem \ref{thm3} that the
efficiency of KTS depends on $\theta$.

Since the bounding polygons $P_1$ and $P_2$ are defined by the coefficients of $f$,
our approach to derive $\theta$ is to first derive $\xi$, a function of $m$ and $n$, satisfying
\[
\norm{c_{ij}} \leq \xi \max_{l \leq u,v \leq h} \norm{f(u,v)},
\]
for any coefficient $c_{ij}$ of $f$.
But the following lemma shows that one needs only derive
the equivalent of $\xi$ for univariate polynomial to derive $\xi$ itself.

\begin{lemma}
\label{uni_theta_to_bi}
Assume there exists a function $h(n)$ such that
\begin{equation}
\label{asuni}
\norm{b_{i}} \leq h(n) \max_{l \leq t \leq h} \norm{g(t)}
\end{equation}
for any $b_i$ $(i = 0,1,\ldots,n)$, and any univariate polynomial $g(t) = \sum_{i=0}^n b_i \phi_i(t)$.
Then
\begin{equation}
\norm{c_{ij}} \leq h(m)h(n) \max_{l \leq u,v \leq h} \norm{f(u,v)},
\end{equation}
for any $c_{ij}$ $(i = 0,1,\ldots,m;$ $j=0,1,\ldots,n)$, and any
bivariate polynomial $f(u,v) = \sum_{i=0}^m \sum_{j=0}^n c_{ij} \phi_i(u)\phi_j(v)$.
\begin{proof}


Let $f(u,v) = \sum_{i=0}^m \sum_{j=0}^n c_{ij} \phi_i(u)\phi_j(v)$ be an arbitrary bivariate polynomial.  For any $i_0 = 0, 1, \ldots, m$, define $g_{i_0}(v) = \sum_{j=0}^n c_{i_0, j} \phi_j(v)$.  Applying (\ref{asuni}) to $g_{i_0}(v)$ yields
\begin{equation}
\label{e1}
\norm{c_{i_0,j}} \leq h(n)\max_{l \leq v \leq h} \norm{g_{i_0}(v)},
\end{equation}
for any $j = 0, 1, \ldots, n$.  Let $v^*_{i_0} = \arg\max_{l \leq v \leq h} \norm{g_{i_0}(v)}$.  Define $l_{i_0}(u) = \sum_{i=0}^m \sum_{j=0}^n c_{ij} \phi_i(u) \phi_j(v^*_{i_0})$.  Applying (\ref{asuni}) to $l_{i_0}(v)$ yields
\begin{equation}
\label{e2}
\norm{\sum_{j=0}^n c_{ij}\phi_j(v^*_{i_0})} \leq h(m)\max_{l \leq u \leq h} \norm{l_{i_0}(u)},
\end{equation}
for any $i = 0, 1, \ldots, m$.  Consequently, by combining (\ref{e1}) and (\ref{e2}), \\ $\norm{c_{i_0,j}} \leq h(n)\max_{l \leq v \leq h} \norm{g_{i_0}(v)} = h(n) \norm{g_{i_0}(v^*_{i_0})} = h(n)\norm{\sum_{j=0}^n c_{i_0,j}\phi_j(v^*_{i_0})} \leq h(m)h(n)\max_{l \leq u \leq h} \norm{l_{i_0}(u)} \leq h(m)h(n)\max_{l \leq u,v \leq h} \norm{f(u,v)}$.
\end{proof}
\end{lemma}

We proceed to derive $\xi$ of the three bases, starting with the Bernstein basis.
The following lemma regarding the product
of two polynomials in Bernstein basis is needed to find $\xi$
for Bernstein case.

\begin{lemma}
\label{productbez} Let
\[
\begin{array}{llll}
f(t) &=& \sum_{i=0}^n c_{i} Z_{i,n}(t),
& 0 \leq t \leq 1
\end{array}
\]
and
\[
\begin{array}{llll}
g(t) &=& \sum_{i=0}^{n'} c'_{i}
Z_{i,n'}(t), & 0 \leq t \leq 1.
\end{array}
\]
Then
\[
f(t)g(t) = \sum_{i=0}^{n+n'} b_{i}
Z_{i,n+n'}(t),
\]
where
\[
|b_{i}| \leq \max_{i} |c_{i}| \cdot \max_{i} |c'_{i}|.
\]
\begin{proof}
Straightforward arithmetic shows that
\[
b_{i} = \sum_{k=\max(0,i-n')}^{\min(n,i)}
         \frac{\nchoosek{n}{k}\nchoosek{n'}{i-k}}{\nchoosek{n+n'}{i}}
         c_{k}c'_{i-k}.
\]
Taking absolute value on both sides and bounding $|c_{k}|$ (resp.
$|c'_{i-k}|$) with $\max_{i} |c_{i}|$ (resp. $\max_{i}
|c'_{i}|$) gives
\[
|b_{i}| \leq \max |c_{i}| \cdot \max |c'_{i}|
\sum_{k=\max(0,i-n')}^{\min(n,i)}
         \frac{\nchoosek{n}{k}\nchoosek{n'}{i-k}}{\nchoosek{n+n'}{i}}.
\]
Recall the combinatorial identity
\[
\nchoosek{n+n'}{i} =
\sum_{k=\max(0,i-n')}^{\min(n,i)}\nchoosek{n}{k}\nchoosek{n'}{i-k}.
\]
Hence, the lemma follows.
\end{proof}
\end{lemma}

With the above lemma, we are ready to derive $\xi$ of the Bernstein basis.

\begin{thm}
\label{bibound} Let $f(t)$ be a polynomial system
\[
\begin{array}{llll}
f(t) &=& \sum_{i=0}^n c_{i} Z_{i,n}(t), &
0 \leq t \leq 1,
\end{array}
\]
where $c_{i} \in \mathbb{R}^d$. The norm of the coefficients
can be bounded by
\begin{equation}
\label{ceq}  \norm{c_{i}} \leq  \xi_B(n) \max_{t : 0
\leq t \leq 1} \norm{f(t)},
\end{equation}
where
\[
\xi_B(n) = \sum_{i=0}^n \prod_{j=0,1,\ldots,i-1,i+1,\ldots,n} \frac{\max\{|n-j|,|j|\}}{|i-j|}
= O(n^{n+1}).
\]
\end{thm}

\begin{flushleft}
\textbf{Remark.} An inequality in the other
direction, namely, that
\[
\max_{t : 0 \leq t \leq 1} \norm{f(t)} \leq \max
\norm{c_{i}},
\]
is a well-known consequence of the convex hull property of
Bernstein polynomials \cite{farin}.
\end{flushleft}

\begin{proof}
By definition of infinity norm, it suffices to prove the lemma for
the case $c_{i} \in \mathbb{R}$.  Therefore, it is assumed that
$d = 1$ for the rest of this proof.

Let $t_j = j/n$ $(j = 0,1,\ldots,n)$. Define a matrix $A \in \mathbb{R}^{(n+1)
\times (n+1)}$ having element
\[
A_{j+1,i+1}=Z_{i,n}(t_j).
\]
Define the vectors $c = (c_0, c_1, \ldots, c_n)^T$ and
$f = \left( f(t_0), f(t_1), \ldots, f(t_n) \right)^T.$
Observe that
\begin{equation}
\label{abf} Ac = f.
\end{equation}
We claim that $A$ is invertible. In particular, we show that the
linear system $Ax = b$ has solution for any arbitrary $b \in
\mathbb{R}^{n+1}$. Due to the definition of $A$, solving
the system $Ax = b$ is equivalent to finding the coefficients
of the polynomial
\begin{equation}
\label{1dbez} g(t) = \sum_{i=0}^n x_{i+1} Z_{i,n}(t)
\end{equation}
with the property that $g(t_0) = b_1$, $g(t_1) = b_2$,
\ldots, $g(t_n) = b_{n+1}$.  The polynomial $g$ satisfying such
property is the Lagrange interpolant
\begin{equation}
\label{lagint} g(t) = \sum_{j=0}^n \left( b_{j+1} \prod_{j'
= 0,\ldots,j-1,j+1,\ldots,n} \frac{t-t_{j'}}{t_j-t_{j'}}
                          \right).
\end{equation}
Transforming (\ref{lagint}) to the Bernstein basis
yields the solution $x$.

Knowing that $A$ is invertible, we multiply both sides of
(\ref{abf}) by $\inv{A}$,
\begin{equation}
c = \inv{A}f,
\end{equation}
and hence, for any $i = 0, 1, \ldots, n$,
\begin{eqnarray}
|c_{i}| & \leq & \norm{c} \nonumber \\
                       & \leq & \norm{\inv{A}}\cdot\norm{f} \nonumber \\
                       & \leq & \norm{\inv{A}}\cdot\max_{t : 0 \leq t \leq 1} \abs{f(t)}.
\label{beforec}
\end{eqnarray}
Comparing (\ref{ceq}) to (\ref{beforec}), it is seen that the
final step is to show that $\norm{\inv{A}} \leq \xi_B(n)$.

Observe that the $i$th column of $\inv{A}$ is $\inv{A} e_i$, where
$e_i$ denotes the $i$th column of the identity matrix. Let
$g_i(t)$ be a polynomial in the Bernstein basis and let $\{ c'_{i'} \}$ be its coefficients. With
similar reasoning as the above,
\begin{equation}
\label{bbb} \left( \begin{array}{c}
      c'_{0} \\ \vdots \\ c'_{n} \end{array} \right)
= \inv{A} \left( \begin{array}{c}
      g_i(t_0) \\ \vdots \\ g_i(t_n) \end{array} \right).
\end{equation}
But (\ref{bbb}) implies that the $i$th column of $\inv{A}$,
$\inv{A}e_i$, are the coefficients of $g_i$ such that, for
$j=0,1,\ldots, n$,
\begin{equation}
\label{gcon} g_i(t_{j}) = \left\{ \begin{array}{ll}
                     1, & j = i,    \\
                     0, & j \neq i. \end{array}  \right.
\end{equation}
The following Lagrange interpolant $g_i$
satisfies (\ref{gcon}):
\begin{eqnarray}
g_i(t) & = & \prod_{j = 0,\ldots,i-1,i+1,\ldots,n}
\frac{t-t_{j}}{t_i-t_{j}} \nonumber \\
& = & \prod_{j = 0,\ldots,i-1,i+1,\ldots,n}
\left( \frac{n-j}{i-j}t - \frac{j}{i-j}(1-t) \right). \label{lastg}
\end{eqnarray}
Note that each term of the product in (\ref{lastg}) is a polynomial in
Bernstein basis with coefficients  $(n-j)/(i-j)$ and $j/(i-j)$.
Applying Lemma \ref{productbez} to (\ref{lastg}) shows
that
\begin{equation}
\label{bp} \norm{\inv{A}e_i} \leq  \prod_{j=0,1,\ldots,i-1,i+1,\ldots,n} \frac{\max\{|n-j|,|j|\}}{|i-j|}.
\end{equation}
Since (\ref{bp}) holds for any column $i$ of $\inv{A}$,
the lemma follows.
\end{proof}

Next is the derivation of $\xi$ of the Chebyshev basis.  The following
identity is useful for this derivation:
\begin{equation}
\sum_{k=1}^n T_i(t_k)T_j(t_k) = \left\{\begin{array}{ll}0 & i \neq j \\ n & i=j=0\\
n/2 & i=j\neq 0,
\end{array} \right.
\label{relation_zeros_cheb}
\end{equation}
for $i,j=0,\ldots,n-1$, where $t_k$ $(k=1,2,\ldots,n)$ are the $n$ zeros of $T_n(t)$.

\begin{thm}
\label{chebbound} Let $f(t)$ be a polynomial system
\[
f(t) = \sum_{i=0}^n c_{i} T_i(t),
\]
where $c_{i} \in \mathbb{R}^d$. The norm of the coefficients
can be bounded by
\begin{equation}
\norm{c_{i}} \leq  \sqrt{2} \max_{t : -1
\leq t \leq 1} \norm{f(t)}.
\end{equation}
\end{thm}
\begin{proof}
By definition of infinity norm, it suffices to prove the lemma for
the case $c_i \in \mathbb{R}$.  Therefore, it is assumed that
$d = 1$ for the rest of this proof.

Let $t_j$ $(j=1,2,\ldots,n+1)$ be the $n+1$ zeros of $T_{n+1}(t)$,
which lie in $[-1,1]$. Define a matrix $A \in \mathbb{R}^{(n+1)
\times (n+1)}$ having element
\[
A_{j+1,i+1}=T_i(t_j).
\]
Define the vectors $c = (c_0, c_1, \ldots, c_n)^T$ and
$f = \left( f(t_0), f(t_1), \ldots, f(t_n) \right)^T.$
Observe that
\begin{equation}
\label{abf_cheb}
Ac = f.
\end{equation}
By (\ref{relation_zeros_cheb}),
\[
A^T A = \mathrm{diag}\left(n+1,(n+1)/2,(n+1)/2,\ldots,(n+1)/2\right),
\]
which implies that $A$ is invertible and
\begin{equation}
\label{aat}
\inv{A}A^{-T} = \mathrm{diag}\left(1/(n+1),2/(n+1),2/(n+1),\ldots,2/(n+1)\right).
\end{equation}
The equation (\ref{aat}) implies
\begin{equation}
\normtwo{\inv{A}}=\sqrt{2/(n+1)}.
\label{normAi}
\end{equation}
Finally, from (\ref{abf_cheb}) and (\ref{normAi}),
\begin{eqnarray}
\abs{c_i} & \leq & \normtwo{c} \nonumber \\
& \leq & \normtwo{\inv{A}}\normtwo{f} \nonumber \\
& \leq & \sqrt{n+1}\normtwo{\inv{A}}\norm{f} \nonumber \\
& \leq & \sqrt{n+1}\normtwo{\inv{A}}\max_{t : -1 \leq t \leq 1} \abs{f(t)} \nonumber \\
& = & \sqrt{2}\max_{t : -1 \leq t \leq 1} \abs{f(t)}. \nonumber \qedhere
\end{eqnarray}
\end{proof}

Last is the power basis.  Our approach to derive $\xi$ of power basis
is to derive the relationship between the coefficients of
a polynomial in power basis and the coefficients of the same
polynomial but written in Chebyshev basis.  By using this relationship
and Theorem \ref{chebbound}, $\xi$ of the power basis can be computed.

\begin{lemma}
\label{power_cheb_relation}
Let $f$ be a univariate polynomial such that
\[
f(t) = \sum_{i=0}^n a_i t^i = \sum_{i=0}^n c_i T_i(t).
\]
In other words, $\{ a_i \}$ are the coefficients of $f$ when written in the power
basis and $\{
c_i \}$ are the coefficients of $f$ when written in the Chebyshev basis.
Then
\[
|a_i| \leq \frac{3^{n+1}-1}{2} \max_{j=0,\ldots,n} |c_j|,
\]
for any $i=0,\ldots,n$.
\begin{proof}
Let $D = [d_{i,j}]$ be the $n+1$-by-$n+1$ matrix such that
\[
a = Dc,
\]
where $a = (a_0, a_1, \ldots, a_n)^T$ and $c =
(c_0, c_1, \ldots, c_n)^T$.  Note that
\[
T_j(t) = \sum_{i=0}^j d_{i+1,j+1}t^i.
\]
Recall the recurrence relation $T_j(t) = 2tT_{j-1}(t)-T_{j-2}(t)$.
It follows from this recurrence that
\begin{equation}
|d_{i+1,j+1}| \leq 3^j. \label{d1}
\end{equation}
That is, when $T_j(t)$ is written in the power basis, the resulting
coefficients (of power basis) is less than or equal to $3^j$.  The
inequality (\ref{d1}) can be verified by induction on the recurrence
relation.  Since the entries in the $(j+1)$th column of $D$ is
bounded by $3^j$, we have, from geometric sum,
\[
\norm{D} \leq (3^{n+1}-1)/2.
\]
The lemma follows from $\norm{a} \leq \norm{D}\norm{c}$.
\end{proof}
\end{lemma}
\begin{thm}
\label{powerbound} Let $f(t)$ be a polynomial system
\[
f(t) = \sum_{i=0}^n c_{i}t^i,
\]
where $c_{i} \in \mathbb{R}^d$. The norm of the coefficients
can be bounded by
\begin{equation}
\norm{c_{i}} \leq  \frac{3^{n+1}-1}{\sqrt{2}} \max_{t : -1
\leq t \leq 1} \norm{f(t)}.
\end{equation}
\end{thm}
\begin{proof}
Follow directly from Theorem \ref{chebbound} and Lemma \ref{power_cheb_relation}.
\end{proof}

Having $\xi$ for each of the three bases, the values of $\theta$ for
the three bases can now be derived.

\begin{cor}
Let
\[
f(u,v) = \sum_{i=0}^m \sum_{j=0}^n c_{ij}Z_{i,m}(u)Z_{j,n}(v),
\]
where $c_{ij} \in \mathbb{R}^2$ $(i=0,1,\ldots,m; j=0,1,\ldots,n)$.
Let $P_1$ be the convex hull of $\{ c_{ij} \}$. Then, for any $y \in P_1$,
\[
\norm{y} \leq \xi_B(m) \xi_B(n) \max_{0 \leq u,v \leq 1}\norm{f(u,v)}.
\]
\end{cor}
\begin{proof}
By the convex hull property of Bernstein polynomials, $\norm{y} \leq \max_{i,j} \norm{c_{ij}}$
for any $y \in P_1$.
The corollary then follows from Theorem \ref{bibound} and Lemma \ref{uni_theta_to_bi}.
\end{proof}

\begin{cor}
Let
\[
f(u,v) = \sum_{i=0}^m \sum_{j=0}^n c_{ij}u^i v^j,
\]
where $c_{ij} \in \mathbb{R}^2$ $(i=0,1,\ldots,m; j=0,1,\ldots,n)$.
Let
\[
P_2 = \left\{ c_{00} + \sum_{i+j > 0} c_{ij} s_{ij} : -1
\leq s_{ij} \leq 1 \right\}.
\]
Then, for any $y \in P_2$,
\[
\norm{y} \leq  \frac{(m+1)(n+1)(3^{m+1}-1)(3^{n+1}-1)}{2}\max_{-1 \leq u,v \leq 1}\norm{f(u,v)}.
\]
\end{cor}
\begin{proof}
For any $y \in P_2$,
\[
\norm{y} \leq \sum_{i=0}^m \sum_{j=0}^n \norm{c_{ij}}.
\]
The corollary then follows from Theorem \ref{powerbound} and Lemma \ref{uni_theta_to_bi}.
\end{proof}

\begin{cor}
Let
\[
f(u,v) = \sum_{i=0}^m \sum_{j=0}^n c_{ij}T_i(u)T_j(v),
\]
where $c_{ij} \in \mathbb{R}^2$ $(i=0,1,\ldots,m; j=0,1,\ldots,n)$.
Let
\[
P_2 = \left\{ c_{00} + \sum_{i+j > 0} c_{ij} s_{ij} : -1
\leq s_{ij} \leq 1 \right\}.
\]
Then, for any $y \in P_2$,
\[
\norm{y} \leq  2(m+1)(n+1)\max_{-1 \leq u,v \leq 1}\norm{f(u,v)}.
\]
\end{cor}
\begin{proof}
For any $y \in P_2$,
\[
\norm{y} \leq \sum_{i=0}^m \sum_{j=0}^n \norm{c_{ij}}.
\]
The corollary then follows from Theorem \ref{chebbound} and Lemma \ref{uni_theta_to_bi}.
\end{proof}

\section{Relationship between the bounding polygon of the power basis and that of Chebyshev basis}
\label{section_cheb_better_than_power} Let $P_2^p$ denote the
bounding polygon $P_2$ computed from the power basis
representation of a polynomial and $P_2^c$ denote $P_2$ computed
from the Chebyshev basis representation of it. The results from
previous section show that the value $\theta$ of $P_2^c$ is
smaller than $\theta$ of $P_2^p$. This only implies that the worst
case of $P_2^c$ is better than the worst case of $P_2^p$.
Comparing the values of $\theta$'s of the two does not indicate
that $P_2^c$ is always a better choice than $P_2^p$ for every
polynomial.  The following results show, however, that $P_2^c$ is,
in fact, always a better choice than $P_2^p$. Specifically, this
section shows that for any given polynomial, its bounding polygon
$P_2^c$ is a subset of its bounding polygon $P_2^p$.

The following two lemmas show that when representing monomials
$t^k$ in Chebyshev basis, each coefficient is nonnegative, and the
sum of all coefficients are exactly 1.  These results are useful
in relating $P_2^p$ to $P_2^c$.

\begin{lemma}
\label{lem_nonnegu} Let $d_{ki}$'s ($k = 0,1,\ldots$; $i =
0,1,\ldots, k$) be the numbers satisfying $t^k = \sum_{i=0}^k
d_{ki}T_i(t)$.  Then
\[
d_{ki} \geq 0,
\]
for any $k = 0,1,\ldots$ and any $i = 0,1,\ldots, k$.
\begin{proof}
We prove the lemma by induction on $k$.  The base cases $k=0$ and
$k=1$ are trivial.  For the inductive step, for any $k \geq 1$,
\begin{eqnarray}
t^{k+1} & = & t \cdot t^k \nonumber \\
& = & t \sum_{i=0}^k d_{ki} T_i(t) \nonumber \\
& = & \sum_{i=1}^k \frac{d_{ki}}{2} \left(2tT_i(t)-T_{i-1}(t)
\right)+ \sum_{i=1}^k \frac{d_{ki}}{2}T_{i-1}(t) + d_{k0}tT_0(t).
\nonumber
\end{eqnarray}
By (\ref{cheb_recurrence}) and noting that $tT_0(t) = t = T_1(t)$,
\begin{eqnarray}
t^{k+1} & = & \sum_{i=1}^k \frac{d_{ki}}{2}T_{i+1}(t)+
\sum_{i=1}^k
\frac{d_{ki}}{2}T_{i-1}(t) + d_{k0}T_1(t)  \nonumber \\
& = & \sum_{i=k}^{k+1}\frac{d_{k,i-1}}{2}T_{i}(t) +
\sum_{i=2}^{k-1} \frac{d_{k,i-1}+d_{k,i+1}}{2}T_i(t) +
\left(\frac{d_{k2}}{2}+d_{k0}\right)T_1(t) +
\frac{d_{k1}}{2}T_0(t). \nonumber
\end{eqnarray}
Hence,
\begin{equation}
d_{k+1,i} = \left\{\begin{array}{ll} d_{k,i-1}/2, &
i=k,k+1, \\
\left(d_{k,i-1}+d_{k,i+1}\right)/2, & i=2,\ldots,k-1, \\
d_{k2}/2+d_{k0},& i=1, \\
d_{k1}/2,& i = 0. \end{array} \right. \label{ukplus1def}
\end{equation}
But since $d_{ki} \geq 0$ for any $i=0,\ldots,k$ by the induction
hypothesis, (\ref{ukplus1def}) shows that $d_{k+1,i} \geq 0$ for
any $i=0,\ldots,k+1$.
\end{proof}
\end{lemma}

\begin{lemma}
\label{lem_sumeq1}
Let $d_{ki}$'s ($k = 0,1,\ldots$; $i =
0,1,\ldots, k$) be the numbers satisfying $t^k = \sum_{i=0}^k
d_{ki}T_i(t)$. Then
\[
\sum_{i=0}^kd_{ki} = 1,
\]
for any $k = 0,1,\ldots$ and any $i = 0,1,\ldots, k$.
\begin{proof}
We prove the lemma by induction on $k$.  The base cases $k=0$ and
$k=1$ are trivial.  For the inductive step, the same reasoning as
in the proof of Lemma \ref{lem_nonnegu} shows that $d_{k+1,i}$ is
as in (\ref{ukplus1def}) for any $k \geq 1$.  Therefore,
\begin{eqnarray}
\sum_{i=0}^{k+1}d_{k+1,i} & = &
\sum_{i=k}^{k+1}\frac{d_{k,i-1}}{2} + \sum_{i=2}^{k-1}
\frac{d_{k,i-1}+d_{k,i+1}}{2} +
\left(\frac{d_{k2}}{2}+d_{k0}\right) + \frac{d_{k1}}{2} \nonumber
\\
& = & \sum_{i=0}^k d_{ki} = 1, \nonumber
\end{eqnarray}
by the induction hypothesis.
\end{proof}
\end{lemma}

\begin{thm}
Let $f : \mathbb{R}^2 \rightarrow \mathbb{R}^2$ be a bivariate
polynomial. Its bounding polygon $P_2^c$ is a subset of its
bounding polygon $P_2^p$.
\begin{proof}
Let $f = f(u,v) = \sum_{k=0}^m \sum_{l=0}^n a_{kl}u^kv^l =
\sum_{i=0}^m \sum_{j=0}^n c_{ij}T_i(u)T_j(v)$, where $a_{kl} \in
\mathbb{R}^n$ ($k=0,1,\ldots,m$; $l=0,1,\ldots,n$) are the
coefficients of $f$ when written in the power basis and $c_{ij}
\in \mathbb{R}^n$ ($i=0,1,\ldots,m$; $j=0,1,\ldots,n$) are the
coefficients of $f$ when written in the Chebyshev basis. Let
$d_{ki}$'s ($k = 0,1,\ldots$; $i = 0,1,\ldots, k$) be the numbers
satisfying $t^k = \sum_{i=0}^k d_{ki}T_i(t)$.  Hence,
\begin{eqnarray}
f(u,v) & = & \sum_{k=0}^m \sum_{l=0}^n a_{kl} \left(\sum_{i=0}^k
d_{ki}T_i(u)\right)\left( \sum_{j=0}^l d_{lj}T_j(v) \right)
\nonumber \\
& = & \sum_{i=0}^m \sum_{j=0}^n \sum_{k=i}^m \sum_{l=j}^n  a_{kl}
d_{ki} d_{lj} T_i(u) T_j(v). \nonumber
\end{eqnarray}
Therefore, $c_{ij} = \sum_{k=i}^m \sum_{l=j}^n  a_{kl} d_{ki}
d_{lj}$.  This means that $P_2^c$ can be written as
\[
P_2^c = \left\{a_{00}d_{00}d_{00} + \sum_{k=1}^m \sum_{l=1}^n
a_{kl} d_{k0} d_{l0} + \sum_{i+j > 0} \sum_{k=i}^m \sum_{l=j}^n
a_{kl} d_{ki} d_{lj}s_{ij} : -1 \leq s_{ij} \leq 1\right\},
\]
But since $d_{00} = 1$,
\begin{eqnarray}
P_2^c & = & \left\{a_{00} + \sum_{k=1}^m \sum_{l=1}^n a_{kl}
d_{k0} d_{l0} + \sum_{j=1}^n \sum_{k=0}^m \sum_{l=j}^n a_{kl}
d_{k0} d_{lj}s_{0j} + \right. \nonumber \\
& & \hspace{11 pt}\left. \sum_{i=1}^m \sum_{k=i}^m \sum_{l=0}^n
a_{kl} d_{ki} d_{l0}s_{i0} + \sum_{i=1}^m\sum_{j=1}^n \sum_{k=i}^m
\sum_{l=j}^n a_{kl} d_{ki} d_{lj}s_{ij} : -1 \leq s_{ij} \leq 1\right\} \nonumber \\
& = & \left\{a_{00} + \sum_{k=1}^m \sum_{l=1}^n a_{kl} d_{k0}
d_{l0} + \sum_{k=0}^m \sum_{l=1}^n \sum_{j=1}^l  a_{kl}
d_{k0} d_{lj}s_{0j} + \right. \nonumber \\
& & \hspace{11 pt}\left. \sum_{k=1}^m \sum_{l=0}^n \sum_{i=1}^k
a_{kl} d_{ki} d_{l0}s_{i0} + \sum_{k=1}^m
\sum_{l=1}^n \sum_{i=1}^k\sum_{j=1}^l  a_{kl} d_{ki} d_{lj}s_{ij} : -1 \leq s_{ij} \leq 1\right\} \nonumber \\
& = & \left\{ a_{00} + \sum_{l=1}^n  a_{0l} \sum_{j=1}^l
d_{lj}s_{0j} + \sum_{k=1}^m a_{k0} \sum_{i=1}^k d_{ki}
s_{i0} + \right. \nonumber \\
& & \hspace{11 pt}\left. \sum_{k=1}^m \sum_{l=1}^n a_{kl} \left(
d_{k0}d_{l0} + d_{k0} \sum_{j=1}^l d_{lj}s_{0j} +
d_{l0}\sum_{i=1}^k d_{ki} s_{i0} + \sum_{i=1}^k d_{ki}\sum_{j=1}^l
d_{lj}s_{ij}
\right) : \right. \nonumber \\
& & \hspace{11 pt}-1 \leq s_{ij} \leq 1 \Bigg\}. \nonumber
\end{eqnarray}
By Lemma \ref{lem_nonnegu} and Lemma \ref{lem_sumeq1}, it is seen
that $-1 \leq \sum_{j=1}^l d_{lj}s_{0j}, \sum_{i=1}^k d_{ki}
s_{i0} \leq 1$.  In addition, using the fact that $\abs{s_{ij}}
\leq 1$, for any $i=0,\ldots,m$ and any $j=0,\ldots,n$, together
with Lemma \ref{lem_nonnegu} and Lemma \ref{lem_sumeq1}, it is
seen that \\$\abs{d_{k0}d_{l0} + d_{k0} \sum_{j=1}^l d_{lj}s_{0j}
+ d_{l0}\sum_{i=1}^k d_{ki} s_{i0} + \sum_{i=1}^k
d_{ki}\sum_{j=1}^l d_{lj}s_{ij}} \leq \abs{d_{k0}d_{l0}} + $ \\
$\abs{d_{k0}} \sum_{j=1}^l \abs{d_{lj}} + \abs{d_{l0}}\sum_{i=1}^k
\abs{d_{ki}} + \sum_{i=1}^k \abs{d_{ki}}\sum_{j=1}^l \abs{d_{lj}}
= d_{k0}d_{l0} + d_{k0} \sum_{j=1}^l d_{lj} + d_{l0}\sum_{i=1}^k
d_{ki} + \sum_{i=1}^k d_{ki}\sum_{j=1}^l d_{lj} =
\left(\sum_{i=0}^k d_{ki} \right)\left(\sum_{j=0}^l d_{lj} \right)
= 1$.  Therefore, $P_2^c \subseteq P_2^p$.
\end{proof}
\end{thm}

\subsection{Reparametrization}
The last nontrivial basis property that warrants detailed discussion
is the issue of efficient reparametrization.
Reparametrizing polynomials
in power basis is straightforward from the \emph{binomial} theorem.
Polynomials in other bases, on
the other hand, may not be as simple to reparametrize.  The
details of the process for polynomials in Bernstein and Chebyshev
bases are covered in this section.

\subsubsection{Reparametrization of polynomials in Bernstein
basis}  There is more than one algorithm to compute the
reparametrization with $[0,1]^2$ of a bivariate polynomial in
Bernstein basis. We describe one method here. Our method makes use
of a program that, given $\alpha_{ij}$'s, $c$, $d$, $e$, $f$, $g$,
$h$, $k$, and $l$, computes $\beta_{ij}$'s satisfying
\[
\sum_{i=0}^m \sum_{j=0}^n
\alpha_{ij}(cy+d)^i(ey+f)^{m-i}(gz+h)^j(kz+l)^{n-j} =  \sum_{i=0}^m \sum_{j=0}^n \beta_{ij}y^iz^j.
\]
Such conversion can be done in $O\left((mn)^2\right)$ by
generalizing Horner's rule. We leave the details of the conversion
to the reader. Let $X$ denote $\{(u,v):u^0 - r \leq u \leq u^0+r,
v^0 - r \leq v \leq v^0+r \}$. To compute the coefficients
$\{\hat{c}_{ij}\}$ of $\{\hat{f}_X(\hat{u},\hat{v}) : 0 \leq
\hat{u}, \hat{v} \leq 1\}$, the $[0,1]^2$-reparametrized surface
of $\{f(u,v): u^0 - r \leq u \leq u^0+r, v^0 - r \leq v \leq v^0+r
\}$, first substitute $u = 2r \hat{u}+u^0-r$ and $v =
2r\hat{v}+v^0-r$ into $f$, yielding
\begin{eqnarray}
f(u,v) & = &
\sum_{i=0}^m \sum_{j=0}^n \nchoosek{m}{i}\nchoosek{n}{j}c_{ij}u^i(1-u)^{m-i}v^j(1-v)^{n-j} \nonumber \\
& = & \sum_{i=0}^m \sum_{j=0}^n \nchoosek{m}{i}\nchoosek{n}{j}c_{ij}(2r \hat{u}+u^0-r)^i\left(1-(2r \hat{u}+u^0-r)\right)^{m-i} \cdot \nonumber \\
& & \hspace{40 pt}(2r \hat{v}+v^0-r)^j\left(1-(2r \hat{v}+v^0-r)\right)^{n-j} \nonumber \\
& = & \sum_{i=0}^m \sum_{j=0}^n \nchoosek{m}{i}\nchoosek{n}{j}c_{ij}\left((u^0+r)\hat{u}+(u^0-r)(1-\hat{u})\right)^i
\cdot \nonumber \\
& & \hspace{40 pt}\left((1-u^0-r)\hat{u}+(1-u^0+r)(1-\hat{u})\right)^{m-i} \cdot \nonumber \\
& & \hspace{40 pt}\left((v^0+r)\hat{v}+(v^0-r)(1-\hat{v})\right)^j \cdot \nonumber \\
& & \hspace{40 pt}\left((1-v^0-r)\hat{v}+(1-v^0+r)(1-\hat{v})\right)^{n-j}. \label{l2}
\end{eqnarray}
Substituting $\hat{u} = \tilde{u}/(\tilde{u}+1)$ and $\hat{v} = \tilde{v}/(\tilde{v}+1)$
into (\ref{l2}) yields
\begin{eqnarray}
f(u,v) &=& \frac{1}{(\tilde{u}+1)^m (\tilde{v}+1)^n} \sum_{i=0}^m \sum_{j=0}^n \nchoosek{m}{i}\nchoosek{n}{j}c_{ij}
\left( (u^0+r)\tilde{u} + u^0 - r  \right)^i \cdot \nonumber \\
& & \hspace{130 pt} \left((1-u^0-r)\tilde{u}+1-u^0+r \right)^{m-i} \cdot \nonumber \\
& & \hspace{130 pt} \left( (v^0+r)\tilde{v} + v^0 - r  \right)^j \cdot \nonumber \\
& & \hspace{130 pt} \left((1-v^0-r)\tilde{v}+1-v^0+r \right)^{n-j}. \label{l3} \\
& = & \frac{1}{(\tilde{u}+1)^m (\tilde{v}+1)^n}\sum_{i=0}^m \sum_{j=0}^n \gamma_{ij}\tilde{u}^i\tilde{v}^j, \label{l4}
\end{eqnarray}
where (\ref{l4}) is obtained from (\ref{l3}) by
the conversion program mentioned above. Substituting $\hat{u}$ and $\hat{v}$ back
into (\ref{l4}) to see that
\begin{equation}
f(u,v) =
\sum_{i=0}^m \sum_{j=0}^n \gamma_{ij}\hat{u}^i(1-\hat{u})^{m-i}\hat{v}^j(1-\hat{v})^{n-j}.
\end{equation}
Therefore, $\hat{c}_{ij} = \gamma_{ij}/(C(m,i)C(n,j))$ are the
control points of $\hat{f}_X$ where $C(m,i) = \nchoosek{m}{i}$.

\subsubsection{Reparametrization of polynomials in Chebyshev
basis}  Let $a$, $b$, $d$ and $e$ be scalar constants.  The
reparametrization with $[-1,1]^2$ of a bivariate polynomial in
Chebyshev basis can be computed if the values of $\lambda_{ik}$'s
$(i = 0,1,\ldots, m)$ satisfying
\[
T_i(at+b) = \sum_{k=0}^i \lambda_{ik} T_k(t)
\]
and the values of $\mu_{jk}$'s $(j=0,1,\ldots,n)$ satisfying
\[
T_j(dt+e) = \sum_{k=0}^j \mu_{jk} T_k(t)
\]
are known.  Note that
\[
T_i(au+b)T_j(dv+e) = \sum_{k=0}^i \sum_{k'=0}^j \lambda_{ik}
\mu_{jk'} T_k(u)T_{k'}(v),
\]
which is adequate to find the reparametrization.  The values of
$a$, $b$, $d$, and $e$ are determined by the $uv$-domain of the
surface to be reparametrized.

To compute $\lambda_{ik}$'s, observe that for $i \geq 1$, by
(\ref{cheb_recurrence}),
\begin{eqnarray}
T_{i+1}(at+b) &=& 2(at+b)T_i(at+b) - T_{i-1}(at+b) \nonumber \\
& = & 2(at+b)\sum_{k=0}^i \lambda_{ik}T_k(t) - \sum_{k=0}^{i-1}
\lambda_{i-1,k}T_k(t) \nonumber \\
& = & \sum_{k=0}^i 2a\lambda_{ik}tT_k(t) + \sum_{k=0}^i 2b
\lambda_{ik}T_k(t) - \sum_{k=0}^{i-1}
\lambda_{i-1,k}T_k(t) \nonumber \\
& = & 2a\lambda_{i0}tT_0(t) + \sum_{k=1}^i a\lambda_{ik}
\left(2tT_k(t)-T_{k-1}(t) \right) + \nonumber \\
& & \sum_{k=0}^{i-1} a\lambda_{i,k+1} T_k(t) + \sum_{k=0}^i 2b
\lambda_{ik}T_k(t) - \sum_{k=0}^{i-1}
\lambda_{i-1,k}T_k(t) \nonumber \\
& = & 2a\lambda_{i0}T_1(t) + \sum_{k=1}^i a \lambda_{ik}
T_{k+1}(t) + \nonumber \\
& & \sum_{k=0}^{i-1} a\lambda_{i,k+1} T_k(t) + \sum_{k=0}^i 2b
\lambda_{ik}T_k(t) - \sum_{k=0}^{i-1} \lambda_{i-1,k}T_k(t),
\label{c1}
\end{eqnarray}
and
\begin{eqnarray}
T_0(at+b) &=& T_0(t),  \label{c2} \\
T_1(at+b) &=& aT_1(t)+bT_0(t). \label{c3}
\end{eqnarray}
The equalities (\ref{c1}), (\ref{c2}), and (\ref{c3}) yield a
recurrence relation of $\lambda_{ik}$'s that can be used to
compute their values.  The values of $\mu_{jk}$'s can be computed
similarly.

\section{Computational results}

Three versions of KTS algorithms are implemented in Matlab; one
operating on the polynomials in power basis, one on Bernstein
basis, and one on Chebyshev basis.  They are tested against a
number of problem instances with varying condition numbers.  Most
of the test problems are created by using the normally distributed
random numbers as the coefficients $c_{ij}$'s of $f$ in Chebyshev
basis.  For some of the test problems especially those with high
condition number, some coefficients are manually entered.  The
resulting Chebyshev polynomial system is then transformed to the
equivalent system in the power and the Bernstein bases. Hence the
three versions of KTS solve the same polynomial system and the
efficiency of the three are compared. The degrees of the test
polynomials are between biquadratic and biquartic.

For the experiment, we use the algorithm by J\'{o}nsson and
Vavasis \cite{jonsson} to compute the complex zeros required to
estimate the condition number. Table \ref{table_ba} compares the
efficiency of the three versions of KTS for the test problems
with differing condition numbers. The total number of subpatches
examined by KTS during the entire computation and the width of
the smallest patch among those examined are reported. The results
do not show any one version to be particularly more efficient than
the others although the Chebyshev basis has better theoretical
bound than the other two.
\begin{table}
\caption{Comparison of the efficiency of KTS algorithm
operating on the power, the Bernstein, and the Chebyshev bases.
The number of patches examined during the course of the algorithm
and the width of the smallest patch examined are shown for each
version of KTS. \label{table_ba} }
\begin{center}
\begin{tabular}{|r|r|r|r|r|r|r|}
\hline
 & \multicolumn{2}{c}{Power basis} \vline& \multicolumn{2}{c}{Bernstein basis}\vline & \multicolumn{2}{c}{Chebyshev basis}\vline\\
\cline{2-7}
$\cond{f}$ & Num. of & Smallest & Num. of & Smallest & Num. of & Smallest\\
& \multicolumn{1}{c}{patches}\vline &
\multicolumn{1}{c}{width}\vline &
\multicolumn{1}{c}{patches}\vline &
\multicolumn{1}{c}{width}\vline
& \multicolumn{1}{c}{patches}\vline & \multicolumn{1}{c}{width}\vline\\
\hline \hline
$3.8 \times 10^3$ & 29 & .125 & 17 & .0625 & 21 & .125 \\
$1.3 \times 10^4$ & 13 & .125 & 17 & .0625 & 13 & .125 \\
$2.5 \times 10^5$ & 49 & .0625 & 21 & .0625 & 45 & .0625 \\
$1.1 \times 10^6$ & 97 & .0313 & 65 & .0313 & 85 & .0313 \\
$3.9 \times 10^7$ & 89 & .0313 & 81 & .0313 & 89 & .0313 \\
\hline
\end{tabular}
\end{center}
\end{table}

Since the types of test polynomials may affect the relative
efficiency of the three versions of KTS, another experiment is
performed on degree $6$ univariate polynomials generated by
different methods.  Since Section
\ref{section_cheb_better_than_power} shows that the Chebyshev
basis always gives tighter bounding polygons than the power basis,
this experiment only compares between the Chebyshev and Bernstein
bases.  Table \ref{table_bb} and Table \ref{table_bc} show the
results of this experiment.  The polynomials are generated as
follows. The ``rand" polynomials are generated by interpolating
points whose x-coordinates are evenly spaced between $-1$ and $1$,
inclusive, and whose y-coordinates are normally distributed random
numbers. The ``sin" ones are interpolations of $\sin(ax+b)$ at
evenly spaced points between $-1$ and $1$, inclusive, where $a$
and $b$ are normally distributed random numbers.  The ``sin-L"
ones are the same as the ``sin" ones except that least-squares
interpolation is used instead. The ``sinw" (resp. ``sinw-L") ones
are generated in the same way as the ``sin" (resp. ``sin-L") ones
but with the function $\sin(6ax+b)$. Table \ref{table_bb} compares
the number of test polynomials of each type where one basis yields
tighter bounding intervals than the other.  Table \ref{table_bc}
shows the number of test polynomials of each type that bounding
intervals of each basis have at least one endpoint exactly at the
boundary of the ranges of the polynomials. The results show that
the Chebyshev basis is decidedly better for ``rand" polynomials,
is about the same for ``sin" ones, but is worse for the rests of
the polynomials than the Bernstein basis.
\begin{table}
\caption{The numbers of test polynomials out of $1000$ that
bounding intervals associated with the Bernstein basis is tighter
than the those associated with the Chebyshev basis, and \emph{vice
versa}. \label{table_bb} } \begin{center}
\begin{tabular}{|r|r|r|}
\hline Poly. type & Num. that Bernstein is tighter & Num. that Chebyshev is tighter\\
\hline \hline
rand & 1 & 999 \\
sin & 963 & 37 \\
sin-L & 960 & 40 \\
sinw & 436 & 564 \\
sinw-L & 998 & 2 \\
\hline
\end{tabular}
\end{center}
\end{table}

\begin{table}
\caption{The numbers of test polynomials out of $1000$ that
bounding intervals associated with the Bernstein basis and those
associated with the Chebyshev basis having at least one endpoint
exactly at the boundary of the ranges of the polynomials.
\label{table_bc} }
\begin{center}
\begin{tabular}{|r|r|r|}
\hline Poly. type & Num. Bernstein with & Num. Chebyshev with \\
& \multicolumn{1}{c}{ exact endpoint}\vline & \multicolumn{1}{c}{ exact endpoint}\vline  \\
 \hline \hline
rand & 2 & 13 \\
sin & 965 & 0 \\
sin-L & 972 & 0\\
sinw & 330 & 0 \\
sinw-L & 999 & 658 \\
\hline
\end{tabular}
\end{center}
\end{table}

\section{Conclusion}

Three common bases, the power, the Bernstein, and the Chebyshev
bases, are shown to satisfy the required properties for KTS to
perform efficiently. In particular, the values of $\theta$ for the
three bases are derived. These values are used to calculate the
time complexity of KTS when that basis is used to represent the
polynomial system.  The Chebyshev basis has the smallest $\theta$
among the three, which shows that using KTS with the Chebyshev
basis has the smallest worst-case time complexity.  The
computational results, however, show no significant differences
between the performances of the three versions of KTS operating
on the three bases. It appears that, in average case, choosing any
of the three bases do not greatly affect the efficiency of KTS.
The experiment on univariate polynomials show that the Bernstein
basis is more suitable for certain types of polynomials while the
Chebyshev basis is better suited for other types.

\bibliography{basis}

\end{document}